\documentclass[12pt,leqno]{article}
\usepackage{amssymb,graphics,latexsym}
\usepackage{graphicx}
\usepackage{amsmath}

\def\R{\hbox{\bf\rlap{I}{\hbox to 2 pt{}}R}}

\newcommand{\re}{{\rm Re\, }}
\newcommand{\im}{{\rm Im\, }}

\newcommand{\dia}{{\rm diag\, }}

\begin{document}
\thispagestyle{empty}
\begin{center}
{\LARGE \bf Zero-dilation Index of a Finite Matrix}
\end{center}

\vspace*{3mm}

\centerline{\bf Hwa-Long Gau$^{\mbox{a, *, 1}}$, \hspace*{.3cm} \bf Kuo-Zhong Wang$^{\mbox{b, 1}}$,\hspace*{.3cm}\bf Pei Yuan Wu$^{\mbox{b, 1}}$}

\vspace*{3mm}

\noindent
${ }^{\mbox{a}}$Department of Mathematics, National Central University, Chungli 32001, Taiwan \\
${ }^{\mbox{b}}$Department of Applied Mathematics, National Chiao
Tung University, Hsinchu \\ \hspace*{1.5mm} 30010, Taiwan

\vspace{3mm}

\noindent
{\bf Abstract}

For an $n$-by-$n$ complex matrix $A$, we define its zero-dilation index $d(A)$ as the largest size of a zero matrix which can be dilated to $A$. This is the same as the maximum $k$ ($\ge 1$) for which 0 is in the rank-$k$ numerical range of $A$. Using a result of Li and Sze, we show that if $d(A) > \lfloor 2n/3\rfloor$, then, under unitary similarity, $A$ has the zero matrix of size $3d(A)-2n$ as a direct summand. It complements the known fact that if $d(A)>\lfloor n/2\rfloor$, then 0 is an eigenvalue of $A$. We then use it to give a complete characterization of $n$-by-$n$ matrices $A$ with $d(A)=n-1$, namely, $A$ satisfies this condition if and only if it is unitarily similar to $B\oplus 0_{n-3}$, where $B$ is a 3-by-3 matrix whose numerical range $W(B)$ is an elliptic disc and whose eigenvalue other than the two foci of $\partial W(B)$ is 0. We also determine the value of $d(A)$ for any normal matrix and any weighted permutation matrix $A$.

\vspace*{3mm}

\hspace*{19mm} permutation matrix

\vspace*{3mm}

${}^*$Corresponding author.

\emph{E-mail addresses}: hlgau@math.ncu.edu.tw (H.-L. Gau), kzwang@math.nctu.edu.tw (K.-Z. Wang), pywu@math.nctu.edu.tw (P. Y. Wu).

${ }^{1}$Research supported by the National Science Council of the Republic of China under the NSC-101-2115-M-008-006,  NSC-101-2115-M-009-001 and NSC-101-2115-M-009-004 projects, respectively. The third author also acknowledges the support from the MOE-ATU project.

\newpage

\noindent
{\large\bf 1. Introduction}

\vspace{5mm}

Let $A$ be an $n$-by-$n$ complex matrix. In this paper, we define the \emph{zero-dilation index} of $A$ by
$$d(A)=\max\{k\ge 1 : 0_k \mbox{ dilates to } A\},$$
where $0_k$ denotes the $k$-by-$k$ zero matrix. Recall that a $k$-by-$k$ matrix $B$ is said to \emph{dilate} to $A$ (or $B$ is a \emph{compression} of $A$) if $B=V^*AV$ for some $n$-by-$k$ matrix $V$ with $V^*V=I_k$, the $k$-by-$k$ identity matrix, or, equivalently, if $A$ is unitarily similar to a matrix of the form $\left[{\scriptsize\begin{array}{cc} B & \ast\\ \ast & \ast\end{array}}\right]$. Another way to express $d(A)$ is via the totally isotropic subspaces of $A$. Note that a subspace $M$ of $\mathbb{C}^n$ is \emph{totally isotropic} for $A$ if $\langle Ax, y\rangle=0$ for all $x$ and $y$ in $M$, where $\langle\cdot,\cdot\rangle$ denotes the standard inner product in $\mathbb{C}^n$. Thus $d(A)$ is the same as the maximum dimension of the totally isotropic subspaces of $A$.

\vspace{5mm}

The notion of the zero-dilation index is closely related to that of the higher-rank numerical range. Recall that for an $n$-by-$n$ matrix $A$ and a $k$, $1\le k\le n$, the \emph{rank-$k$ numerical range} $\Lambda_k(A)$ of $A$ is the subset $\{\lambda\in\mathbb{C} : \lambda I_k \mbox{ dilates to } A\}$ of the complex plane. In particular, $\Lambda_1(A)$ is simply the classical \emph{numerical range} $W(A)$ of $A$. The study of the higher-rank numerical ranges is a hotly pursued one in recent years starting with the paper \cite{2} by Choi, Kribs and \.{Z}yczkowski. It is known that the $\Lambda_k(A)$'s are always convex (cf. \cite[Theorem 0.1]{12} or \cite[Corollary 2.3]{10}). Obviously, $d(A)$ is equal to the maximum $k$ for which $\Lambda_k(A)$ contains $0$. In proving the convexity of the $\Lambda_k(A)$'s, Li and Sze gave a more specific description of $\Lambda_k(A)$ \cite[Theorem 2.2]{10}, namely,
$$\Lambda_k(A)=\bigcap_{\theta\in\mathbb{R}}\{\lambda\in\mathbb{C}: \re(e^{-i\theta}\lambda)\le\lambda_k(\re(e^{-i\theta}A))\},$$
where, for a complex number $z$ and a matrix $B$, $\re z=(z+\overline{z})/2$ and $\re B=(B+B^*)/2$ are their \emph{real parts}, and, for an $n$-by-$n$ Hermitian matrix $C$, $\lambda_1(C)\ge \cdots\ge \lambda_n(C)$ denote its eigenvalues arranged in decreasing order. In terms of this description, they also gave in \cite[Theorem 3.1]{10} the expression
$$d(A)=\min\{k_{\theta} : \lambda_{k_{\theta}}(\re(e^{-i\theta}A))\ge 0,  \lambda_{k_{\theta}+1}(\re(e^{-i\theta}A))< 0, \theta\in\mathbb{R}\}$$
for $d(A)$.

\vspace{5mm}

In Section 2 below, we first give some basic properties of the zero-dilation index, some of which are based on the Li--Sze theorem. For example, we show in Proposition 2.1 that if $A$ is an $n$-by-$n$ matrix with 0 in $\partial W(A)$, then $d(A)\le\dim\bigvee\{x\in\mathbb{C}^n: \langle Ax, x\rangle=0\}$ and, moreover, the equality holds if and only if 0 is an extreme point of $W(A)$. Theorem 2.2 is the Li--Sze theorem and Corollary 2.4 is an easy consequence of it. The latter says that $d(A)=\min\{d(\re(e^{-i\theta}A)) : \theta\in\mathbb{R}\}$ for any matrix $A$, which essentially reduces the computation of the zero-dilation index of a general matrix to those of Hermitian matrices. The latter can be done quite easily as in Corollary 2.3.

\vspace{5mm}

Section 3 relates large values of $d(A)$ to the zero eigenvalue of $A$. In particular, Theorem 3.2 says that if $A$ is of size $n$ and $d(A)>\lfloor 2n/3\rfloor$, then $A$ is unitarily similar to a matrix of the form $B\oplus 0_{3d(A)-2n}$ and the number $\lfloor 2n/3\rfloor$ is sharp. This is in contrast to the situation for $d(A)>\lfloor n/2\rfloor$, in which case we only have 0 as an eigenvalue of $A$ (cf. \cite[Proposition 2.2]{2}). Using the former, we characterize in Theorem 3.3 those $n$-by-$n$ matrices $A$ with $d(A)=n-1$: this is the case if and only if $A$ is unitarily similar to $B\oplus 0_{n-3}$, where $B$ is a 3-by-3 matrix whose numerical range $W(B)$ is an elliptic disc and whose eigenvalues are 0 and the two foci of the ellipse $\partial W(B)$.

\vspace{5mm}

Finally, in Section 4, we determine the zero-dilation indices of normal and weighted permutation matrices. If $A$ is an $n$-by-$n$ normal matrix and $k$ is the number of nonzero eigenvalues of $A$, then $d(A)$ is an integer between $k$ and $\lfloor (n+k)/2\rfloor$. We also characterize those normal $A$'s with their $d(A)$'s attaining the extremal $k$ and $\lfloor (n+k)/2\rfloor$ (cf. Theorem 4.1). Since every weighted permutation matrix (a square matrix with at most one nonzero entry on each of its rows and columns) is permutationally similar to the direct sum of matrices of the forms
$$\left[\begin{array}{cccc} 0 & w_1 & & \\ & 0 & \ddots & \\ & & \ddots & w_{n-1}\\ & & & 0\end{array}\right] \ \ \ \mbox{and} \ \ \  \left[\begin{array}{cccc} 0 & w_1 & & \\ & 0 & \ddots & \\   & & \ddots & w_{n-1}\\ w_n & & & 0\end{array}\right],$$
where the $w_j$'s are all nonzero, its zero-dilation index can be determined from those of the latter two types. It turns out that
$$d(A)=\left\{\begin{array}{ll} \lceil\frac{1}{2}n\rceil & \mbox{      if } A \mbox{ is of the first type},\\ \lfloor\frac{1}{2}n\rfloor & \mbox{      if } A \mbox{ is of the second type}\end{array}\right.$$
(cf. Corollary 2.7 and Lemma 4.4, respectively), and hence the zero-dilation index of a weighted permutation matrix can be computed as in Theorem 4.5.

\vspace{5mm}

We end this section by fixing some notations. For any $m$-by-$n$ matrix $A$, $A^T$ (resp., $A^*$) denotes its transpose (resp., adjoint). We use $0_{m\, n}$ to denote the $m$-by-$n$ zero matrix; this is abbreviated to $0_n$ if $m=n$. If $A$ is a square matrix, then $\re A=(A+A^*)/2$ and $\im A=(A-A^*)/(2i)$ are its \emph{real} and \emph{imaginary parts}, respectively. Two $n$-by-$n$ matrices $A$ and $B$ are \emph{permutationally similar} if there is a \emph{permutation matrix} $V$, that is, one with exactly one 1 on each of its rows and columns, such that $V^*AV=B$. We use $\dia(a_1, \ldots, a_n)$ to denote the $n$-by-$n$ diagonal matrix with the diagonals $a_1, \ldots, a_n$. For any subset $K$ of $\mathbb{C}^n$,  $\bigvee K$ denotes the subspace of $\mathbb{C}^n$ generated by the vectors in $K$. If $z$ is a nonzero complex number, then $\theta\equiv\arg z$ is the unique number in $[0, 2\pi)$ satisfying $z=|z|e^{i\theta}$. If $x$ is a real number, then $\lfloor x\rfloor$ (resp., $\lceil x\rceil$) denotes the largest (resp., smallest) integer less than (resp., greater than) or equal to $x$. For any set $\bigtriangleup$, $\#\bigtriangleup$ denotes its cardinality. If $\bigtriangleup$ is a Lebesgue measurable subset of $\mathbb{R}$, then $|\bigtriangleup|$ denotes its Lebesgue measure. For a subset $\bigtriangleup$ of $\mathbb{C}$, $\bigtriangleup^{\wedge}$ denotes its convex hull.

\vspace{5mm}

Our reference for general properties of numerical ranges of matrices is \cite[Chapter 1]{6}.

\vspace{8mm}

\noindent
{\large\bf 2. Preliminaries}

\vspace{5mm}

We start with the following proposition for the value of $d(A)$ when 0 is in the boundary of the numerical range of $A$.

\vspace{5mm}

{\bf Proposition 2.1.} \emph{If $A$ is an $n$-by-$n$ matrix with $0$ in $\partial W(A)$},
\emph{then $d(A)\le\dim\bigvee\{x\in \mathbb{C}^n : \langle Ax, x\rangle=0\}$}. \emph{Moreover}, \emph{in this case}, \emph{the equality holds if and only if $0$ is an extreme point of $W(A)$}.

\vspace{5mm}

Recall that a point $\lambda$ is an \emph{extreme point} of the convex subset $\bigtriangleup$ of the plane if $\lambda$ is in $\bigtriangleup$ and it cannot be expressed as $t\lambda_1+(1-t)\lambda_2$ with $\lambda_1$ and $\lambda_2$ in $\bigtriangleup$ both distinct from $\lambda$ and $0<t<1$.

\vspace{5mm}

{\em Proof of Proposition $2.1$}.  Let $d=d(A)$, $K=\{x\in \mathbb{C}^n : \langle Ax, x\rangle=0\}$ and $k=\dim \bigvee K$. Since $U^*AU=\left[{\scriptsize\begin{array}{cc} 0_d & \ast\\ \ast & \ast\end{array}}\right]$ for some $n$-by-$n$ unitary matrix $U$, we have $\langle Ax, x\rangle=0$ for all $x$ in $M\equiv U(\mathbb{C}^d\oplus\{0\})$. This shows that $M\subseteq K\subseteq\bigvee K$ and hence $d=\dim M\le\dim\bigvee K=k$.

\vspace{3mm}

If $d=k$, then $M=K=\bigvee K$ from above. Hence $K$ is a subspace of $\mathbb{C}^n$, which is equivalent to 0 being an extreme point of $W(A)$ (cf. \cite[Theorem 1 (i)]{3}). Conversely, if 0 is extreme for $W(A)$, then $K$ is a subspace of $\mathbb{C}^n$.  The compression $A_1=P_{K}A|K : K\rightarrow K$ of $A$ to $K$, where $P_K$ is the (orthogonal) projection of $C^n$ onto $K$, is such that $\langle A_1x, x\rangle=\langle P_KAx, x\rangle=\langle Ax, x\rangle =0$ for all $x$ in $K$. Hence we deduce that $A_1=0_k$ and, therefore, $A$ is unitarily similar to  $\left[{\scriptsize\begin{array}{cc} 0_k & \ast\\ \ast & \ast\end{array}}\right]$. The maximality of $d$ implies that $k\le d$. Together with the already-proven $d\le k$, this yields their equality.  $\blacksquare$

\vspace{5mm}

Next we reformulate \cite[Theorem 3.1]{10} in terms of our terminology. For a Hermitian matrix $A$, let $i_+(A)$ (resp., $i_-(A)$ and $i_0(A)$) denote the number of positive (resp., negative and zero) eigenvalues of $A$ (counting multiplicity), $i_{\ge 0}(A)=i_+(A)+i_0(A)$, and $i_{\le 0}(A)=i_-(A)+i_0(A)$.

\vspace{5mm}

{\bf Theorem 2.2.} \emph{For any $n$-by-$n$ matrix $A$}, \emph{we have} $d(A)=\min\{i_{\ge 0}(\re(e^{-i\theta}A)):\theta\in\mathbb{R}\}.$

\vspace{5mm}

Several corollaries follow, some of which are inspired by the results in \cite{2} and \cite{10}.

\vspace{5mm}

{\bf Corollary 2.3.} \emph{If $A$ is an $n$-by-$n$ Hermitian matrix}, \emph{then} $d(A)=\min\{i_{\ge 0}(A), i_{\le 0}(A)\}$. \emph{In particular}, \emph{in this case}, $i_0(A)\le d(A)\le\lfloor (n+i_0(A))/2\rfloor$.

\vspace{5mm}

{\em Proof}. It is obvious that
$$i_+(\re(e^{-i\theta}A))=\left\{\begin{array}{ll} i_+(A) & \mbox{      if } 0\le\theta<\frac{\pi}{2} \mbox{  or  } \frac{3\pi}{2}<\theta<2\pi,\\
0 & \mbox{      if } \theta=\frac{\pi}{2} \mbox{   or   } \frac{3\pi}{2},\\ i_-(A) & \mbox{      if } \frac{\pi}{2}<\theta<\frac{3\pi}{2}, \end{array} \right.$$
and
$$i_0(\re(e^{-i\theta}A))=\left\{\begin{array}{ll} i_0(A) & \mbox{      if } 0\le\theta<2\pi \mbox{  and  } \theta\neq \frac{\pi}{2},  \frac{3\pi}{2},\\ n & \mbox{      if } \theta=\frac{\pi}{2} \mbox{ or } \frac{3\pi}{2}. \end{array} \right.$$
Our assertions on $d(A)$ then follow immediately from Theorem 2.2.   $\blacksquare$

\vspace{5mm}

The preceding bounds on $d(A)$ will be extended to a normal matrix $A$ in Theorem 4.1.

\vspace{5mm}

{\bf Corollary 2.4.} \emph{For any $n$-by-$n$ matrix $A$}, \emph{we have} $d(A)=\min\{d(\re(e^{-i\theta}A)):\theta\in\mathbb{R}\}.$

\vspace{5mm}

{\bf Corollary 2.5.} \emph{If $A$ is an $n$-by-$n$ matrix such that $\dim\ker(\re(e^{-i\theta_0}A))\le 1$ for some real $\theta_0$}, \emph{then} $d(A)\le\lceil n/2\rceil$.

\vspace{5mm}

{\em Proof}. Let $k=i_{\ge 0}(\re(e^{-i\theta_0}A))$. If $k\le\lceil n/2\rceil$, then Theorem 2.2 implies that $d(A)\le k\le\lceil n/2\rceil$. On the other hand, if $k>\lceil n/2\rceil$, then, since $i_-(\re(e^{-i\theta_0}A))<n-\lceil n/2\rceil$ and $i_0(\re(e^{-i\theta_0}A))\le 1$, we have
\begin{eqnarray*}
d(A)&\le&i_{\ge 0}(\re(e^{-i(\theta_0+\pi)}A))=i_{\le 0}(\re(e^{-i\theta_0}A))\\
&<& (n-\lceil\frac{1}{2}n\rceil)+1 \le \lceil\frac{1}{2}n\rceil+1
\end{eqnarray*}
by Theorem 2.2. Hence, in this case, $d(A)\le\lceil n/2\rceil)$ holds too. $\blacksquare$

\vspace{5mm}

Note that the converse of Corollary 2.5 is false as witness the matrix $A=\dia(0,0,1)$ with $d(A)=2$ and $\dim\ker(\re(e^{-i\theta}A))\ge 2$ for all real $\theta$. However, in one case, it is indeed true.

\vspace{5mm}

{\bf Corollary 2.6.} \emph{If $A$ is an $n$-by-$n$ matrix which is unitarily similar to $-A$}, \emph{then} $d(A)\ge \lceil n/2\rceil$. \emph{In this case}, $d(A)=\lceil n/2\rceil$ \emph{if and only if} $\dim\ker\re(e^{-i\theta_0}A)\le 1$ \emph{for some real} $\theta_0$.

\vspace{5mm}

{\em Proof}. Our assumption implies that $\re(e^{-i\theta}A)$ and $-\re(e^{-i\theta}A)$ are unitarily similar and thus $d(\re(e^{-i\theta}A))\ge \lceil n/2\rceil$ for all real $\theta$. Using Corollary 2.4, we obtain $d(A)\ge\lceil n/2\rceil$.

\vspace{3mm}

If $d(A)=\lceil n/2\rceil$, then, by Corollary 2.4, there is some real $\theta_0$ such that $d(\re(e^{-i\theta_0}A))=\lceil n/2\rceil$. The unitary similarity of $\re(e^{-i\theta_0}A)$ and $-\re(e^{-i\theta_0}A)$ yields that $\dim\ker\re(e^{-i\theta_0}A)\le 1$. The converse follows from Corollary 2.5 and what was proven in the preceding paragraph.  $\blacksquare$

\vspace{5mm}

The next corollary gives a class of matrices which satisfy the conditions in Corollary 2.6.

\vspace{5mm}

{\bf Corollary 2.7.}  \emph{If} $$A=\left[\begin{array}{cccc} 0 & w_1 & & \\ & 0 & \ddots & \\ & & \ddots & w_{n-1}\\ & & & 0\end{array}\right]$$ \emph{is of size} $n$ $(\ge 1)$ \emph{with $w_j\neq 0$ for all $j$}, \emph{then} $d(A)=i_{\ge 0}(\re(e^{-i\theta}A))=\lceil n/2\rceil$ \emph{for all real} $\theta$.

\vspace{5mm}

{\em Proof}. It is easily seen that $A$ is unitarily similar to $e^{-i\theta}A$ for all real $\theta$ and $\dim\ker\re A\le 1$. Thus $i_{\ge 0}(\re(e^{-i\theta}A))$ is independent of the value of $\theta$ and hence is equal to $d(A)$ for all $\theta$ by Theorem 2.2. Applying Corollary 2.6, we obtain $d(A)=\lceil n/2\rceil$.  $\blacksquare$

\vspace{5mm}

The zero-dilation indices of general weighted permutation matrices will be determined in Section 4.

\vspace{5mm}

We end this section with the following elementary observation on the zero-dilation index of a direct sum.

\vspace{5mm}

{\bf Corollary 2.8.} \emph{Let} $A=\sum_{j=1}^m\oplus A_j$, \emph{where} $A_j$, $1\le j\le m$, \emph{is an $n_j$-by-$n_j$ matrix}. \emph{Then} $d(A)\ge \sum_{j=1}^m d(A_j)$, \emph{and $d(A)=\sum_{j=1}^m d(A_j)$ if and only if there is some real $\theta_0$ such that $d(A_j)=i_{\ge 0}(\re(e^{-i\theta_0}A_j))$ for all $j$}. \emph{In particular}, $d(\sum_{j=1}^m\oplus A_{1})=m d(A_{1})$ \emph{and} $d(A_1\oplus 0_n)=d(A_1)+n$.

\vspace{5mm}

{\em Proof}. That $d(A)\ge\sum_j d(A_j)$ follows immediately from the definition of the zero-dilation index. Next we assume that $d(A)=\sum_j d(A_j)$. Let the real $\theta_0$ be such that $d(A)=i_{\ge 0}(\re(e^{-i\theta_0}A))$. This is the same as $\sum_jd(A_j)=\sum_j i_{\ge 0}(\re(e^{-i\theta_0}A_j))$. Since $d(A_j)\le i_{\ge 0}(\re(e^{-i\theta_0}A_j))$ for all $j$ by Theorem 2.2, we infer from above that $d(A_j)=i_{\ge 0}(\re(e^{-i\theta_0}A_j))$ for all $j$. For the converse, we need only show that $d(A)\le\sum_j d(A_j)$. For the given $\theta_0$, we have
\begin{eqnarray*}
d(A)&\le& i_{\ge 0}(\re(e^{-i\theta_0}A))=\sum_ji_{\ge 0}(\re(e^{-i\theta_0}A_j))\\
&=& \sum_j d(A_j)
\end{eqnarray*}
as required.

\vspace{3mm}

The assertions for $d(\sum_{j=1}^m\oplus A_1)$ and $d(A_1\oplus 0_n)$ follow easily from above.  $\blacksquare$

\vspace{8mm}

\noindent
{\large\bf 3. Zero Eigenvalue}

\vspace{5mm}

In this section, we consider the relations between large values of $d(A)$ and the zero eigenvalue of $A$. This we start with the following known fact from \cite[Proposition 2.2]{2}.

\vspace{5mm}

{\bf Lemma 3.1.} \emph{If $A$ is an $n$-by-$n$ matrix with} $d(A)>\lfloor n/2\rfloor$, \emph{then $0$ is an eigenvalue of $A$ with} (\emph{geometric}) \emph{multiplicity at least $2d(A)-n$}. \emph{Moreover}, \emph{in this case}, \emph{the number} ``$\lfloor n/2\rfloor$'' \emph{is sharp}.

\vspace{5mm}

{\em Proof}. We only need to show the sharpness of $\lfloor n/2\rfloor$. This is seen by the $n$-by-$n$ matrix
$$A=\left[\begin{array}{cccc} 0 & 1 & & \\ & 0 & \ddots & \\ & & \ddots & 1\\ 1 & & & 0\end{array}\right].$$
Since $A$ is unitarily similar to $\dia(1, \omega_n, \omega_n^2, \ldots, \omega_n^{n-1})$, where $\omega_n=e^{2\pi i/n}$, we infer that, for odd $n$ (resp., even $n$) and for $\theta$ in $[\pi/2, 5\pi/2)$,
$$d(\re(e^{-i\theta}A))=\left\{\begin{array}{ll} \lceil \frac{1}{2}n\rceil \ \ (\mbox{resp., } (\frac{1}{2}n)+1) & \mbox{  if } \theta=\frac{\pi}{2}+j\frac{2\pi}{n}, 0\le j\le n-1,\\ \lfloor\frac{1}{2}n\rfloor \ \ (\mbox{resp., } \frac{1}{2}n) & \mbox{  otherwise}.\end{array}\right.$$
Hence $d(A)=\lfloor n/2\rfloor$ by Corollary 2.4.  $\blacksquare$

\vspace{5mm}

The next theorem says that an even larger value of $d(A)$ will guarantee that 0 be a \emph{reducing eigenvalue} of $A$, meaning that $Ax=A^*x=0$ for some nonzero vector $x$.

\vspace{5mm}

{\bf Theorem 3.2.} \emph{If $A$ is an $n$-by-$n$ matrix with} $d(A)>\lfloor 2n/3\rfloor$, \emph{then $0$ is a reducing eigenvalue of $A$ with} (\emph{geometric}) \emph{multiplicity at least} $3d(A)-2n$ \emph{and $A$ is unitarily similar to a matrix of the form} $B\oplus 0_{3d(A)-2n}$, \emph{where $B$ is of size $3(n-d(A))$ with $d(B)=2(n-d(A))$}. \emph{In this case}, \emph{the number} ``$\lfloor 2n/3\rfloor$'' \emph{is sharp}.

\vspace{5mm}

{\em Proof}. Let $d=d(A)$, $A=[a_{ij}]_{i,j=1}^n$, where $a_{ij}=0$ for $1\le i, j\le d$, and $u_j=[\overline{a}_{d+j, 1} \ \ldots \ \overline{a}_{d+j, d}]^T$ and $v_j=[a_{1, d+j} \ \ldots \ a_{d, d+j}]^T$ for $1\le j\le n-d$. If $W_1=\bigvee\{u_1, \ldots, u_{n-d}\}$ and $W_2=\bigvee\{v_1, \ldots, v_{n-d}\}$, then $W_1$ and $W_2$ are subspaces of $\mathbb{C}^d$, whose orthogonal complements $W_1^{\perp}$ and $W_2^{\perp}$ in $\mathbb{C}^d$ satisfy
\begin{eqnarray*}
\dim(W_1^{\perp}\cap W_2^{\perp})&=& \dim W_1^{\perp} +\dim W_2^{\perp} - \dim(W_1^{\perp} +W_2^{\perp})\\
&\ge& (d-(n-d))+(d-(n-d))-d\\
&=& 3d-2n.
\end{eqnarray*}
Let $x_1, \ldots, x_{3d-2n}$ be orthonormal vectors in $W_1^{\perp}\cap W_2^{\perp}$ and let
$$y_j=\left[\begin{array}{c}x_j \\ 0 \\  \vdots \\  0\end{array}\right]\hspace{-7.2mm}\begin{array}{c} \\ \left. \begin{array}{l}  \\   \\  \end{array}\right\}n-d\end{array}$$
for $1\le j\le 3d-2n$. Then the $y_j$'s are orthonormal in $\mathbb{C}^n$ and $Ay_j=A^*y_j=0$ for all $j$. This yields our assertion on 0 being a reducing eigenvalue of $A$. Hence $A$ is unitarily similar to a matrix of the form $B\oplus 0_{3d(A)-2n}$. That $d(B)=2(n-d(A))$ is a consequence of Corollary 2.8.

\vspace{3mm}

The sharpness of $\lfloor 2n/3\rfloor$ is seen by the $n$-by-$n$ matrix $A=[a_{ij}]_{i,j=1}^n$ with $a_{ij}=1$ if $(i, j)=(n-k+1, k+1), (n-k+2, k+2), \ldots, (k, 3k-n), (k+1, 1), (k+2, 2), \ldots, (n, n-k)$, and $a_{ij}=0$ otherwise, where $k=\lfloor 2n/3\rfloor$. Note that $\ker A$ (resp., $\ker A^*$) consists of vectors of the form $[\underbrace{0 \ \ldots \ 0}_{n-k} \ \underbrace{\ast \ \ldots \ \ast}_{2k-n} \ \underbrace{0 \ \ldots \ 0}_{2k-n} \ \underbrace{\ast \ \ldots \ \ast}_{2n-3k}]^T$ (resp., $[\underbrace{\ast \ \ldots \ \ast}_{n-k} \  \underbrace{0 \ \ldots \ 0}_{k}]^T$). Hence $\ker A \cap \ker A^*=\{0\}$ and, therefore, 0 is not a reducing eigenvalue of $A$. >From what was proven in the first paragraph, we deduce that $d(A)\le \lfloor 2n/3\rfloor=k$. On the other hand, we also have $d(A)\ge k$ by our construction of $A$ and the definition of $d(A)$. Thus $d(A)=\lfloor 2n/3\rfloor$ as required.  $\blacksquare$

\vspace{5mm}

Using the preceding theorem, we can now give a characterization of $n$-by-$n$ matrices $A$ with $d(A)=n-1$.

\vspace{5mm}

{\bf Theorem 3.3.} \emph{Let $A$ be an $n$-by-$n$ $(n\ge 3)$ matrix}. \emph{Then $d(A)=n-1$ if and only if $A$ is unitarily similar to a matrix of the form $B\oplus 0_{n-3}$}, \emph{where $B$ is of size $3$ whose numerical range $W(B)$ is an elliptic disc} (\emph{or a line segment}) \emph{and whose eigenvalues are $0$ and the two foci} (\emph{or the two endpoints}) \emph{of} $\partial W(B)$.

\vspace{5mm}

The next lemma is a special case of Theorem 3.3 for $n=3$.

\vspace{5mm}

{\bf Lemma 3.4.} \emph{Let $A$ be a $3$-by-$3$ matrix}. \emph{Then $d(A)=2$ if and only if $W(A)$ is an elliptic disc} (\emph{or a line segment}) \emph{and the eigenvalues of $A$ are $0$ and the two foci} (\emph{or the two endpoints}) \emph{of} $\partial W(A)$.

\vspace{5mm}

For its proof, we need the Kippenhahn polynomial of a matrix. If $A$ is an $n$-by-$n$ matrix, then its \emph{Kippenhahn polynomial} is $p_A(x,y,z)=\det(x\re A+y\im A+zI_n)$ for $x, y$ and $z$ in $\mathbb{C}$. Note that $p_A$ codifies the information of the spectrum and numerical range of $A$: the roots of $p_A(1, i, -z)=0$ in $z$ are the eigenvalues of $A$ while the convex hull of the real points of the dual curve of $p_A(x,y,z)=0$ is the numerical range of $A$ (cf. \cite[Theorem 10]{8}).

\vspace{5mm}

{\em Proof of Lemma $3.4$}. Assume first that $d(A)=2$. Then $A$ is unitarily similar to a matrix of the form ${\scriptsize\left[\begin{array}{cc} 0_2 & \ast\\ \ast & \ast\end{array}\right]}$, and thus the same is true for $\re A$ and $\im A$. Hence $p_A(x,y,z)=z q(x,y,z)$ for some real homogeneous polynomial $q$ of degree 2. Therefore, $W(A)$ is the convex hull of the point 0 and the real points of the dual curve of $q(x,y,z)=0$. We denote the convex hull of the latter set by $\bigtriangleup$. Then $\bigtriangleup$ is either an elliptic disc or a line segment depending on whether $q$ is irreducible or otherwise. We claim that 0 must be in $\bigtriangleup$. Indeed, if otherwise, then let $\theta$ in $[0, 2\pi)$ be such that $e^{-i\theta}\bigtriangleup$ is in the open right half-plane. We infer that $i_{\le 0}(\re(e^{-i\theta}A))=1$ and thus $d(A)\le d(\re(e^{-i\theta}A))\le 1$ by Corollaries 2.3 and 2.4, which contradicts our assumption of $d(A)=2$. Hence we have $W(A)=\bigtriangleup$. Since the characteristic polynomial of $A$ is $p_A(-1, -i, z)=z q(-1, -i, z)$, the eigenvalues of $A$ are 0 and the two foci (or the two endpoints) of $\partial W(A)$.

\vspace{3mm}

Conversely, if $A$ satisfies the asserted properties, then $p_A(x,y,z)=z q(x,y,z)$, where $q$ is a real homogeneous polynomial of degree 2 (cf. \cite[Theorem 2.2]{7}). If $x=\cos\theta$ and $y=\sin\theta$ for real $\theta$, then $z$ is a divisor of $p_A(\cos\theta, \sin\theta, z)=\det(\re(e^{-i\theta}A)+zI_3)$. This shows that 0 is an eigenvalue of $\re(e^{-i\theta}A)$ for all real $\theta$. Since 0, being an eigenvalue of $A$, is in the elliptic disc $W(A)$, we infer that $i_{\ge 0}(\re(e^{-i\theta}A))=2$ for all $\theta$. Thus $d(A)=2$ by Theorem 2.2.  $\blacksquare$

\vspace{5mm}

We are now ready to prove Theorem 3.3.

\vspace{5mm}

{\em Proof of Theorem $3.3$}. In view of Lemma 3.4, we need only consider for $n\ge 4$. If $d(A)=n-1$, then $d(A)>\lfloor 2n/3\rfloor$. Theorem 3.2 then yields that $A$ is unitarily similar to a matrix of the form $B\oplus 0_{n-3}$, where $B$ is of size 3 with $d(B)=2$. Lemma 3.4 then furnishes the proof of the asserted necessary condition on $B$. The converse also follows easily from Lemma 3.4 and Corollary 2.8.  $\blacksquare$

\vspace{5mm}

{\bf Corollary 3.5.} \emph{Let $A$ be an $n$-by-$n$ $(n\ge 3)$ nilpotent matrix}. \emph{Then $d(A)=n-1$ if and only if $A$ is unitarily similar to a matrix of the form $\left[{\scriptsize\begin{array}{ccc} 0 & a &  \\ & 0 & b\\ & & 0\end{array}}\right]\oplus 0_{n-3}$ with $a$ and $b$ not both $0$}.

\vspace{5mm}

{\em Proof}. This is an easy consequence of Theorem 3.3, Lemma 3.4 and the fact that a 3-by-3 nilpotent matrix $A$ has its numerical range $W(A)$ equal to a circular disc (centered at the origin) if and only if it is unitarily similar to a nonzero matrix of the form   $\left[{\scriptsize\begin{array}{ccc} 0 & a &  \\ & 0 & b\\ & & 0\end{array}}\right]$ (cf. \cite[Theorem 4.1]{7}).  $\blacksquare$

\vspace{5mm}

The first assertion of the following corollary is due to Linden \cite[Proposition 1]{11} (cf. also \cite[Theorem 2]{1}).

\vspace{5mm}

{\bf Corollary 3.6.} \emph{If $A$ is an $n$-by-$n$ $(n\ge 3)$ matrix with $d(A)=n-1$}, \emph{then $W(A)$ is an elliptic disc} (\emph{or a line segment}). \emph{In this case}, \emph{the number} ``$n-1$'' \emph{is sharp}.

\vspace{5mm}

{\em Proof}. The sharpness of $n-1$ is seen by the matrix $A=B\oplus 0_{n-3}$, where $$B=\left[\begin{array}{ccc} 0 & 1 & 1\\ & 0 & 1\\ & & 0\end{array}\right].$$ In this case, we obviously have $d(B)\ge 1$ and hence $d(A)\ge n-2$. If $d(A)=n-1$, then Theorem 3.3 implies that $W(A)$ is an elliptic disc. On the other hand, it is known that $W(B)$ has a line segment on its boundary and contains 0 in its interior (cf. \cite[Theorem 4.1 (2)]{7}). Thus the same is true for $A$, which contradicts what we have shown above. We conclude that $d(A)=n-2$.   $\blacksquare$

\vspace{8mm}

\noindent
{\large\bf 4. Normal Matrix and Weighted Permutation Matrix}

\vspace{5mm}

In this section, we determine the zero-dilation indices for matrices in two special classes: the normal ones and the weighted permutation ones. We start with the former class.

\vspace{5mm}

{\bf Theorem 4.1.} \emph{If $A$ is an $n$-by-$n$ normal matrix with $k=\dim\ker A$}, \emph{then} $k\le d(A)\le \lfloor (n+k)/2\rfloor$. \emph{Moreover}, \emph{let $\lambda_1, \ldots, \lambda_{n-k}$ be the nonzero eigenvalues of $A$} (\emph{counting multiplicity}) \emph{arranged such that} $\arg\lambda_1\le \cdots\le\arg\lambda_{n-k}$. \emph{Then} (a) $d(A)=k$ \emph{if and only if $0$ is not in the convex hull of} $\{\lambda_1, \ldots, \lambda_{n-k}\}$, \emph{and} (b) $d(A)=\lfloor(n+k)/2\rfloor$ \emph{if and only if}, \emph{for even} $n-k$ (\emph{resp}., \emph{odd} $n-k$), \emph{the condition  $\arg\lambda_{j+((n-k)/2)}-\arg\lambda_j=\pi$ holds for all $j$}, $1\le j\le (n-k)/2$ (\emph{resp}., $\arg\lambda_{j+((n-k-1)/2)}-\arg\lambda_j\le \pi$ \emph{for $1\le j\le (n-k+1)/2$ and $\arg\lambda_{j-((n-k+1)/2)}-\arg\lambda_j\le -\pi$ for $(n-k+3)/2\le j\le n-k$}).

\vspace{5mm}

{\em Proof}. Since $0_k$ is a direct summand of the diagonal form of $A$, we have $d(A)\ge k$. To prove $d(A)\le\lfloor(n+k)/2\rfloor$, let $\theta_0$ in $[0, 2\pi)$ be such that the line $y=x\tan\theta_0$ is not perpendicular to any of the $n-k$ lines connecting the origin and some $\lambda_j$. Then $i_0(\re(e^{-i\theta_0}A))=k$. Hence, by Corollary 2.4,
$$d(A)\le d(\re(e^{-i\theta_0}A)) \le \lfloor\frac{1}{2}(n-k)\rfloor+k = \lfloor\frac{1}{2}(n+k)\rfloor.$$

\vspace{3mm}

To prove (a), note that $d(A)=k$ is equivalent to the existence of a real $\theta_0$ such that $\re(e^{-i\theta_0}\lambda_j)<0$ for all $j$, $1\le j\le n-k$. The latter is easily seen to be the same as $\{\lambda_1, \ldots, \lambda_{n-k}\}^{\wedge}$, the convex hull of  $\{\lambda_1, \ldots, \lambda_{n-k}\}$, not containing 0.

\vspace{3mm}

For the proof of (b), let $m=n-k$ and consider $A$ as $B\oplus 0_k$, where $B=\dia(\lambda_1, \ldots, \lambda_{m})$. In view of Corollary 2.8 and the identity $\lfloor(n-k)/2\rfloor+k=\lfloor(n+k)/2\rfloor$, we need only prove for $B$. Note that, by \cite[Corollary 2.4]{10}, we have
\[\Lambda_{\ell}(B)=\bigcap_{1\le j_1<\cdots<j_{m-\ell+1}\le m}\{\lambda_{j_1}, \ldots, \lambda_{j_{m-\ell+1}}\}^{\wedge}\leqno{(1)}\]
for any $\ell$, $1\le\ell\le m$. First assume that $m$ is even. If $d(B)=m/2$ and if there is some pair $\lambda_j$ and $\lambda_{j+(m/2)}$, $1\le j\le m/2$, with $\arg\lambda_{j+(m/2)}-\arg\lambda_j<\pi$, then 0 is not in $\{\lambda_{j}, \ldots, \lambda_{j+(m/2)}\}^{\wedge}$ and hence not in $\Lambda_{m/2}(B)$ by (1), which contradicts our assumption of $d(B)=m/2$. Similarly, if $\arg\lambda_{j+(m/2)}-\arg\lambda_j>\pi$, then 0 is not in $\{\lambda_{j+(m/2)}, \ldots, \lambda_{m}, \lambda_{1}, \ldots, \lambda_{j}\}^{\wedge}$, which also leads to a contradiction. Hence we have $\arg\lambda_{j+(m/2)}-\arg\lambda_j=\pi$ for all $j$, $1\le j\le m/2$, as required. For the converse, note that any set $\bigtriangleup$ consisting of $(m/2)+1$ many $\lambda_j$'s must contain some pair $\lambda_{j_0}$ and $\lambda_{j_0+(m/2)}$ ($1\le j_0\le m/2$). Hence the assumption of $\arg\lambda_{j_0+(m/2)}-\arg\lambda_{j_0}=\pi$ guarantees that 0 is in $\widehat{\bigtriangleup}$. Thus 0 is in $\Lambda_{m/2}(B)$ by (1), and therefore $d(B)\ge m/2$. Together with the already-proven $d(B)\le m/2$, this yields $d(B)=m/2$.

\vspace{3mm}

Next consider for odd $m$. If $d(B)=\lfloor m/2\rfloor=(m-1)/2$, then 0 is in $\Lambda_{(m-1)/2}(B)$ and hence in $\{\lambda_{j+((m-1)/2)}, \ldots, \lambda_{m}, \lambda_{1}, \ldots, \lambda_{j}\}^{\wedge}$ (resp., $\{\lambda_{j-((m+1)/2)}, \ldots, \lambda_{j}\}^{\wedge}$) for all $j$, $1\le j\le (m+1)/2$ (resp., $(m+3)/2\le j\le m$), by (1). We infer that the asserted argument conditions are satisfied. Conversely, assume that these conditions hold. Let $\bigtriangleup$ be any set consisting of $(m+3)/2$ many $\lambda_j$'s. Then $\bigtriangleup$ must contain some pair $\lambda_{j_0-((m+1)/2)}$ and $\lambda_{j_0}$ ($(m+3)/2\le j_0\le m$). Hence we have $\arg\lambda_{j_0-((m+1)/2)}-\arg\lambda_{j_0}\le -\pi$. On the other hand, $\bigtriangleup$ also contains some $\lambda_{j_1}$, $j_0-((m+1)/2)< j_1<j_0$. From our assumptions, we obtain $\arg\lambda_{j_0}-\arg\lambda_{j_1}\le \pi$ and $\arg\lambda_{j_1}-\arg\lambda_{j_0-((m+1)/2)}\le\pi$. These, together with $\arg\lambda_{j_0}-\arg\lambda_{j_0-((m+1)/2)}\ge \pi$, yield that 0 is in $\{\lambda_{j_0-((m+1)/2)}, \lambda_{j_1}, \lambda_{j_0}\}^{\wedge}$ and hence in $\widehat{\bigtriangleup}$. Thus (1) implies that 0 is in $\Lambda_{(m-1)/2}(B)$. Hence $d(B)\ge (m-1)/2$. Together with the already-proven $d(B)\le (m-1)/2$, this yields $d(B)=(m-1)/2=\lfloor m/2\rfloor$.  $\blacksquare$

\vspace{5mm}

Finally, we consider the zero-dilation indices of weighted permutation matrices. Recall that a \emph{weighted permutation matrix} is one with at most one nonzero entry on each of its
 rows and columns. It is easily seen that every such matrix is permutationally similar to the direct sum of matrices of the forms
 $$\left[\begin{array}{cccc} 0 & w_1 & & \\ & 0 & \ddots & \\ & & \ddots & w_{n-1}\\ & & & 0\end{array}\right] \ \ \ \mbox{and} \ \ \  \left[\begin{array}{cccc} 0 & w_1 & & \\ & 0 & \ddots & \\   & & \ddots & w_{n-1}\\ w_n & & & 0\end{array}\right],$$
where all the $w_j$'s are nonzero. The zero-dilation index of a matrix $A$ of the first type was already given in Corollary 2.7: $d(A)=\lceil n/2\rceil$. The next three lemmas prepare us for the calculation of $d(A)$ when $A$ is of the second type.

\vspace{5mm}

{\bf Lemma 4.2.} \emph{Let} $$A=\left[\begin{array}{cccc} 0 & w_1 & & \\ & 0 & \ddots & \\   & & \ddots & w_{n-1}\\ w_n & & & 0\end{array}\right]$$ \emph{be of size $n$ $(\ge 2)$ with $w_j\neq 0$ for all $j$}, $\alpha=\sum_{j=1}^n\arg w_j$, \emph{and $\lambda_1(\theta)\ge\cdots\ge\lambda_n(\theta)$ be the eigenvalues of $\re(e^{-i\theta}A)$ for each real $\theta$}.

(a) \emph{If $n$ is even}, \emph{then $i_{\ge 0}(\re(e^{-i\theta}A))=(n/2)+1$ or $n/2$ for any real $\theta$}. \emph{In this case}, $i_{\ge 0}(\re(e^{-i\theta}A))=(n/2)+1$\emph{ if and only if} $|w_1w_3\cdots w_{n-1}|=|w_2w_4\cdots w_n|$ \emph{and} $\theta=(\alpha+m\pi)/n$, \emph{where} $m=0, \pm 2, \pm 4, \ldots$ (\emph{resp}., $m=\pm 1, \pm 3, \ldots$) \emph{if $n/2$ is even} (\emph{resp}., $n/2$ \emph{is odd}).

(b) \emph{If $n$ is odd}, \emph{then $i_{\ge 0}(\re(e^{-i\theta}A))=(n+1)/2$ or $(n-1)/2$ for any real $\theta$}. \emph{In this case}, $i_{\ge 0}(\re(e^{-i\theta}A))=(n+1)/2$ \emph{if and only if} $(\alpha-(\pi/2)+2m\pi)/n\le \theta\le(\alpha+(\pi/2)+2m\pi)/n$ (\emph{resp}., $(\alpha+(\pi/2)+2m\pi)/n\le \theta\le(\alpha+(3\pi/2)+2m\pi)/n$) \emph{for some} $m=0, \pm 1, \pm 2, \ldots$ \emph{if $(n-1)/2$ is even} (\emph{resp}., $(n-1)/2$ \emph{is odd}).

\vspace{5mm}

For the proof, we need another lemma.

\vspace{5mm}

{\bf Lemma 4.3.} \emph{If}
$$A=\left[\begin{array}{cccc} 0 & w_1 & & \\ & 0 & \ddots & \\ & & \ddots & w_{n-1}\\ & & & 0\end{array}\right] \ \ \ (\mbox{\em resp}., \   A=\left[\begin{array}{cccc} 0 & w_1 & & \\ & 0 & \ddots & \\   & & \ddots & w_{n-1}\\ w_n & & & 0\end{array}\right])$$
\emph{with $w_j\neq 0$ for all $j$}, \emph{then all the eigenvalues of $\re A$ have multiplicity $1$} (\emph{resp}., \emph{at most $2$}).

\vspace{5mm}

{\em Proof}. Let $\lambda$ be an eigenvalue of $\re A$ and let $x=[x_1 \ \ldots \ x_n]^T$ be such that $(\re A)x=\lambda x$. This yields
$$\frac{1}{2}(w_1x_2+\overline{w}_nx_n)=\lambda x_1,$$ and $$\frac{1}{2}(\overline{w}_{j-1}x_{j-1}+w_{j}x_{j+1})=\lambda x_j, \ \ 2\le j\le n-1.$$
Hence
$$x_2=\frac{2\lambda}{w_1}x_1-\frac{\overline{w}_n}{w_1}x_n\equiv \alpha_2x_1+\beta_2x_n,$$ and
$$x_{j+1}=\frac{2\lambda}{w_j}x_j - \frac{\overline{w}_{j-1}}{w_j}x_{j-1}, \ \ 2\le j\le n-1.$$
The latter yields, by iteration, an expression for $x_{j+1}$, $2\le j\le n-1$, as $x_{j+1}=\alpha_{j+1}x_1+\beta_{j+1}x_n$, where $\alpha_{j+1}$ and $\beta_{j+1}$ are scalars which depend only on $\lambda$ and the $w_j$'s. Let $u=[1 \ \alpha_2 \ \ldots \ \alpha_n]^T$ and $v=[0 \ \beta_2 \ \ldots \ \beta_n]^T$. Then $x$ is a linear combination of $u$ and $v$: $x=x_1u+x_n v$. This shows that the multiplicity of $\lambda$ is at most 2. Moreover, if $w_n=0$, then $\beta_2=\cdots=\beta_n=0$ and hence $x$ is a multiple of $u$. This gives the multiplicity of $\lambda$ as 1.  $\blacksquare$

\vspace{5mm}

{\em Proof of Lemma $4.2$}. (a) Assume that $n$ is even. If $U$ is the $n$-by-$n$ unitary matrix $\dia(1, -1, \ldots, 1, -1)$, then $U^*AU=-A$. It follows that $\re(e^{-i\theta}A)$ is unitarily similar to $-\re(e^{-i\theta}A) $ for any real $\theta$. Thus $\lambda_j(\theta)=-\lambda_{n-j+1}(\theta)$ for $1\le j\le n$. Since the eigenvalues of $\re(e^{-i\theta}A)$ have multiplicity at most 2 by Lemma 4.3, we deduce that $\lambda_j(\theta)>0$ (resp., $\lambda_j(\theta)<0$) for $1\le j\le (n/2)-1$ (resp., $(n/2)+2\le j\le n$). Therefore, $i_{\ge 0}(\re(e^{-i\theta}A))=(n/2)+1$ or $n/2$ depending on whether $\lambda_{n/2}(\theta)=\lambda_{(n/2)+1}(\theta)=0$ or otherwise.

\vspace{3mm}

To determine which value $i_{\ge 0}(\re(e^{-i\theta}A))$ assumes, we make use of the expression of the Kippenhahn polynomial $p_A(x,y,z)$ of $A$ given in \cite[Theorem 4.2]{4} to obtain
\begin{eqnarray*}
&& \det(\re(e^{-i\theta}A))=p_A(\cos\theta, \sin\theta, 0)\\
&=& \frac{1}{2^n}\left[(-1)^{n/2}(|w_1w_3\cdots w_{n-1}|^2+|w_2w_4\cdots w_n|^2)-2|w_1\cdots w_n|\cos(n\theta-\alpha)\right].
\end{eqnarray*}
Note that $\lambda_{n/2}(\theta)=\lambda_{(n/2)+1}(\theta)=0$ if and only if $\det(\re(e^{-i\theta}A))=0$ and, from above, the latter is equivalent to
\begin{eqnarray*}
\cos(n\theta-\alpha)&=& \frac{(-1)^{n/2}}{2}\left(\frac{|w_1w_3\cdots w_{n-1}|}{|w_2w_4\cdots w_n|}+\frac{|w_2w_4\cdots w_n|}{|w_1w_3\cdots w_{n-1}|}\right)\\
&\equiv &  \frac{(-1)^{n/2}}{2}(w+\frac{1}{w}),
\end{eqnarray*}
where $w=|w_1w_3\cdots w_{n-1}|/|w_2w_4\cdots w_n|$. If this equation is to be satisfied, then $(w+(1/w))/2\le 1$, which is the case exactly when $w=1$. Thus we conclude that $i_{\ge 0}(\re(e^{-i\theta}A))=(n/2)+1$ if and only if $|w_1w_3\cdots w_{n-1}|=|w_2w_4\cdots w_n|$ and $\cos(n\theta-\alpha)=(-1)^{n/2}$. The latter condition holds exactly when $\theta$ equals one of the asserted values.

\vspace{3mm}

(b) Now assume that $n$ is odd. Let $B$ be the ($n-1$)-by-($n-1$) matrix
$$\left[\begin{array}{cccc} 0 & w_1 & & \\ & 0 & \ddots & \\ & & \ddots & w_{n-2}\\ & & & 0\end{array}\right].$$
Since $n-1$ is even, $B$ is unitarily similar to $-B$ as in (a). Thus the same is true for $\re(e^{-i\theta}B)$ and $-\re(e^{-i\theta}B)$ for all real $\theta$. Together with Lemma 4.3, this implies that $\lambda_j(\re(e^{-i\theta}B))>0$ (resp., $\lambda_j(\re(e^{-i\theta}B))<0$) for $1\le j\le (n-1)/2$ (resp., $(n+1)/2\le j\le n-1$). Using the interlacing property \cite[Theorem 4.3.8]{5} of the eigenvalues of the $n$-by-$n$ Hermitian matrix $\re(e^{-i\theta}A)$ and its ($n-1$)-by-($n-1$) principal submatrix $\re(e^{-i\theta}B)$, we obtain $\lambda_j(\theta)>0$ (resp., $\lambda_j(\theta)<0$) for $1\le j\le (n-1)/2$ (resp., $(n+3)/2\le j\le n$). Thus, for any real $\theta$, we have $i_{\ge 0}(\re(e^{-i\theta}A))=(n+1)/2$ or $(n-1)/2$ depending on whether $\lambda_{(n+1)/2}(\theta)$ is nonnegative or otherwise.

\vspace{3mm}

Note that $\lambda_{(n+1)/2}(\theta)\ge 0$ if and only if $(-1)^{(n-1)/2}\cos(n\theta-\alpha)\ge 0$. Indeed, as in (a), using the expression of $p_A(x,y,z)$ from \cite[Theorem 4.2]{4}, we have
\begin{eqnarray*}
&& \left(\prod_{j=1}^{(n-1)/2}\lambda_j(\theta)\right)\lambda_{(n+1)/2}(\theta)\left(\prod_{j=(n+3)/2}^n(-\lambda_j(\theta))\right)
=(-1)^{(n-1)/2}\det(\re(e^{-i\theta}A))\\
&=& (-1)^{(n-1)/2}p_A(\cos\theta, \sin\theta, 0)=(-1)^{(n-1)/2}\frac{1}{2^{n-1}}|w_1\cdots w_n|\cos(n\theta-\alpha).
\end{eqnarray*}
Since the first and third products in the first term of the above expression are both (strictly) positive, our assertion follows. We conclude that $i_{\ge 0}(\re(e^{-i\theta}A))=(n+1)/2$ if and only if $(-1)^{(n-1)/2}\cos(n\theta-\alpha)\ge 0$, which is the same as the asserted condition for $\theta$.  $\blacksquare$

\vspace{5mm}

An easy consequence of Lemma 4.2 and Theorem 2.2 is the following.

\vspace{5mm}

{\bf Lemma 4.4.} \emph{If}
$$A=\left[\begin{array}{cccc} 0 & w_1 & & \\ & 0 & \ddots & \\   & & \ddots & w_{n-1}\\ w_n & & & 0\end{array}\right]$$
\emph{is of size $n$ $(\ge 2)$ with $w_j\neq 0$ for all $j$}, \emph{then} $d(A)=\lfloor n/2\rfloor$. \emph{Moreover}, \emph{if $n$ is even} (\emph{resp}., \emph{$n$ is odd}), \emph{then $d(A)=i_{\ge 0}(\re(e^{-i\theta}A))$ for all but finitely many values of $\theta$ in any finite interval of $\mathbb{R}$} (\emph{resp}., \emph{for all $\theta$ in the union of open intervals}
$$\bigcup_{m=-\infty}^{\infty}\left(\frac{1}{n}\Big(\big(\sum_{j=1}^n\arg w_j\big)+\frac{\pi}{2}+2m\pi\Big), \frac{1}{n}\Big(\big(\sum_{j=1}^n\arg w_j\big)+\frac{3\pi}{2}+2m\pi\Big)\right)$$
\emph{or} $$\bigcup_{m=-\infty}^{\infty}\left(\frac{1}{n}\Big(\big(\sum_{j=1}^n\arg w_j\big)-\frac{\pi}{2}+2m\pi\Big), \frac{1}{n}\Big(\big(\sum_{j=1}^n\arg w_j\big)+\frac{\pi}{2}+2m\pi\Big)\right)$$
\emph{depending on whether $(n-1)/2$ is even or odd}).

\vspace{5mm}

We are now ready to compute the zero-dilation index of any weighted permutation matrix.

\vspace{5mm}

{\bf Theorem 4.5.} \emph{Let $A$ be a weighted permutation matrix permutationally similar to a matrix of the form} $(\sum_{j=1}^{p+q}\oplus A_j)\oplus(\sum_{k=1}^rB_k)$, \emph{where} $p, q, r\ge 0$,
$$A_j=\left[\begin{array}{cccc} 0 & a_1^{(j)} & & \\ & 0 & \ddots & \\   & & \ddots & a_{n_j-1}^{(j)}\\ a_{n_j}^{(j)} & & & 0\end{array}\right] \mbox{ \em is of size } n_j \, (\ge 2), \ \ 1\le j\le p+q,$$
\emph{and} $$B_k=\left[\begin{array}{cccc} 0 & b_1^{(k)} & & \\ & 0 & \ddots & \\   & & \ddots & b_{m_k-1}^{(k)}\\   & & & 0\end{array}\right] \mbox{ \em is of size } m_k \, (\ge 1), \ \ 1\le k\le r,$$
\emph{with the weights $a_s^{(j)}$ and $b_t^{(k)}$ all nonzero and the sizes} $n_1, \ldots, n_p$ \emph{odd} (\emph{resp}., $n_{p+1}, \ldots, n_{p+q}$ \emph{even}). \emph{If $\alpha_j=\sum_{s=1}^{n_j}\arg a_s^{(j)}$ for $1\le j\le p+q$}, \emph{then}
$$(2) \hspace{5mm} d(A)=\sum_{j=1}^{p+q}\lfloor \frac{1}{2}n_j\rfloor+\sum_{k=1}^r\lceil\frac{1}{2}m_k\rceil+\min_{\theta\in\mathbb{R}}\#\{j : 1\le j\le p, (-1)^{(n_j-1)/2}\cos(n_j\theta-\alpha_j)>0\}. $$

\vspace{3mm}

{\em Proof}. Note that, for each $k$, $1\le k\le r$, $B_k$ is unitarily similar to $e^{-i\theta}B_k$ for all real $\theta$. Hence the number $i_{\ge 0}(\re(e^{-i\theta}B_k))$ is constant for all the $\theta$'s, and, therefore, for each $k$, $d(B_k)=i_{\ge 0}(\re(e^{-i\theta}B_k))$ for all $\theta$. To prove our assertion, we may assume, in view of Corollary 2.8, that $A=\sum_{j=1}^{p+q}\oplus A_j$.

\vspace{3mm}

From Lemma 4.2, we have, for each real $\theta$,
\begin{eqnarray*}
 i_{\ge 0}(\re(e^{-i\theta}A))&=&\Big[\big(\sum_{j=1}^p\lfloor\frac{1}{2}n_j\rfloor\big)+\#\{j: 1\le j\le p, (-1)^{(n_j-1)/2}\cos(n_j\theta-\alpha_j)\ge 0\}\Big]\\
&& +\Big[\sum_{j=p+1}^{p+q}\big((\frac{1}{2}n_j)-1\big)+\#\{j: p+1\le j\le p+q, \cos(n_j\theta-\alpha_j)\neq 0\}\\
&& +2\#\{j: p+1\le j\le p+q, \cos(n_j\theta-\alpha_j)=0\}\Big].
\end{eqnarray*}
Letting
$$f_j(\theta)=\left\{\begin{array}{ll} (-1)^{(n_j-1)/2}\cos(n_j\theta-\alpha_j) & \mbox{      if } 1\le j\le p,\\
\cos(n_j\theta-\alpha_j) & \mbox{      if } p+1\le j\le p+q,\\
-\cos(n_{j-q}\theta-\alpha_{j-q}) & \mbox{      if } p+q+1\le j\le p+2q,\end{array}\right.$$
and $m(\theta)=\#\{j: 1\le j\le p+2q, f_j(\theta)\ge 0\}$ for real $\theta$, we can express $d(A)$ as
\[ d(A)=\min_{\theta\in\mathbb{R}}i_{\ge 0}(\re(e^{-i\theta}A))=\big(\sum_{j=1}^{p+q}\lfloor\frac{1}{2}n_j\rfloor\big)-q+\min_{\theta\in\mathbb{R}}m(\theta). \leqno{(3)}\]
If $\theta_0\in\mathbb{R}$ is such that $m(\theta_0)=\min_{\theta\in\mathbb{R}}m(\theta)$, we claim that $f_j(\theta_0)\neq 0$ for all $j$, $1\le j\le p+2q$. Indeed, assume that $f_j(\theta_0)>0$ for $1\le j\le p_1$, $f_j(\theta_0)=0$ and $f_j$ is strictly increasing (resp., strictly decreasing) on a neighborhood of $\theta_0$ for $p_1+1\le j\le p_1+p_2$ (resp., $p_1+p_2+1\le j\le p_1+p_2+p_3$), $f_j(\theta_0)<0$ for $p_1+p_2+p_3+1\le j\le p$, $f_j(\theta_0)>0$ for $p+1\le j\le p+q_1$, $f_j(\theta_0)=0$ for $p+q_1+1\le j\le p+q_1+q_2$, and $f_j(\theta_0)<0$ for $p+q_1+q_2+1\le j\le p+q$, where $p_1, p_2, p_3, q_1, q_2\ge 0$ with $p_1+p_2+p_3\le p$ and $q_1+q_2\le q$. Then
\begin{eqnarray*}
m(\theta_0) & = & (p_1+p_2+p_3)+(q_1+q_2)+(q_2+(q-q_1-q_2))\\
& = & p_1+p_2+p_3+q+q_2.
\end{eqnarray*}
Since the $f_j$'s are continuous in $\theta$, there is an $\varepsilon_1>0$ such that $f_j(\theta_0+\varepsilon_1)>0$ for $1\le j\le p_1+p_2$ and $p+1\le j\le p+q_1$, $f_j(\theta_0+\varepsilon_1)<0$ for $p_1+p_2+1\le j\le p$ and $p+q_1+q_2+1\le j\le p+q$, and $f_j(\theta_0+\varepsilon_1)\neq 0$ for $p+q_1+1\le j\le p+q_1+q_2$. Then $m(\theta_0+\varepsilon_1)=p_1+p_2+q$. Since $m(\theta_0)\le m(\theta_0+\varepsilon_1)$ or $p_1+p_2+p_3+q+q_2\le p_1+p_2+q$, we obtain $p_3=q_2=0$. Similarly, there is an $\varepsilon_2>0$ such that $f_j(\theta_0-\varepsilon_2)>0$ for $1\le j\le p_1$ and $p+1\le j\le p+q_1$, and $f_j(\theta_0-\varepsilon_2)<0$ for the remaining $j$'s. Hence $m(\theta_0-\varepsilon_2)=p_1+q$. We infer from $m(\theta_0)\le m(\theta_0-\varepsilon_2)$ that $p_2=0$. This proves our claim. We conclude from above that
$$m(\theta_0)=p_1+q=q+\min_{\theta\in \mathbb{R}}\#\{ j : 1\le j\le p, (-1)^{(n_j-1)/2}\cos(n_j\theta-\alpha_j)>0\} $$
and hence (3) becomes
$$d(A)=\big(\sum_{j=1}^{p+q}\lfloor \frac{1}{2}n_j\rfloor\big)+\min_{\theta\in\mathbb{R}}\#\{j : 1\le j\le p, (-1)^{(n_j-1)/2}\cos(n_j\theta-\alpha_j)>0\} $$
as asserted.   $\blacksquare$

\vspace{5mm}

{\bf Corollary 4.6.} \emph{If $A$ and $B$ are $n$-by-$n$ weighted permutation matrices such that the moduli of their corresponding entries are all equal to each other}, \emph{then} $d(A)=d(B)$.

\vspace{5mm}

{\em Proof}. This is because the expression of $d(A)$ in (2) is independent of the moduli of the entries of $A$.  $\blacksquare$

\vspace{5mm}

{\bf Corollary 4.7.} \emph{Let $A$ be a weighted permutation matrix represented as in Theorem} 4.5, \emph{and let} $d=\sum_{j=1}^{p+q}\lfloor n_j/2\rfloor + \sum_{k=1}^r\lceil m_k/2\rceil$. \emph{Then} $d\le d(A)\le d+\lfloor p/2\rfloor$. \emph{Moreover}, $d(A)=d$ \emph{if and only if} $\cap_{j=1}^pS_j\neq\emptyset$, \emph{where}
\[S_j=\left\{\begin{array}{ll} \bigcup\limits_{m=-\infty}^{\infty}\big(\frac{1}{n_j}(\alpha_j+\frac{\pi}{2}+2m\pi), \frac{1}{n_j}(\alpha_j+\frac{3\pi}{2}+2m\pi)\big) & \mbox{   \em if } \, \frac{1}{2}(n_j-1) \mbox{ \em is even},\\
\bigcup\limits_{m=-\infty}^{\infty}\big(\frac{1}{n_j}(\alpha_j-\frac{\pi}{2}+2m\pi), \frac{1}{n_j}(\alpha_j+\frac{\pi}{2}+2m\pi)\big)& \mbox{   \em if } \, \frac{1}{2}(n_j-1) \mbox{ \em is odd}. \end{array}\right. \leqno{(4)}\]

\vspace{5mm}

{\em Proof}. Assume that $d(A)>d+\lfloor p/2\rfloor$. Then (2) implies that, for any real $\theta$, there are more than $\lfloor p/2\rfloor$ many $j$'s among $1, \ldots, p$ such that $(-1)^{(n_j-1)/2}\cos(n_j\theta-\alpha_j)>0$. Since there are also more than $\lfloor p/2\rfloor$ many $j$'s for which
$$(-1)^{(n_j-1)/2}\cos(n_j\theta-\alpha_j)=-(-1)^{(n_j-1)/2}\cos(n_j(\theta+(\pi/n_j))-\alpha_j)<0.$$
This is certainly impossible. Thus we must have $d(A)\le d+\lfloor p/2\rfloor$.

\vspace{3mm}

Finally, the equivalence condition for $d(A)=d$ follows from Corollaries 2.7 and 2.8 and Lemma 4.4. It is also a consequence of (2) as $\min_{\theta\in\mathbb{R}}\#\{j : 1\le j\le p, (-1)^{(n_j-1)/2}\cos(n_j\theta-\alpha_j)>0\}=0$ means that the minimum 0 is attained at some real $\theta_0$ for which $(-1)^{(n_j-1)/2}\cos(n_j\theta_0-\alpha_j)<0$ for all $j$, $1\le j\le p$ (cf. proof of Theorem 4.5), which is in turn equivalent to $\cap_{j=1}^pS_j\neq\emptyset$.  $\blacksquare$

\vspace{5mm}

Admittedly, for a specific weighted permutation matrix $A$, its $d(A)$ is difficult to compute from the expression (2) in Theorem 4.5. However, at least in two cases, we do have a more precise description of $d(A)$. The first one is for $A$ to have only positive weights.

\vspace{5mm}

{\bf Proposition 4.8.} \emph{Let $A$ be a weighted permutation matrix represented as in Theorem} 4.5. \emph{If all the weights $a_s^{(j)}$ and $b_t^{(k)}$ are} (\emph{strictly}) \emph{positive}, \emph{then}
$$d(A)=\sum_{j=1}^{p+q}d(A_j) + \sum_{k=1}^r d(B_k)=\sum_{j=1}^{p+q}\lfloor\frac{1}{2}n_j\rfloor + \sum_{k=1}^r\lceil \frac{1}{2}m_k\rceil.$$

\vspace{5mm}

{\em Proof}. From our assumption, we have $\alpha_j\equiv\sum_{s=1}^{n_j}\arg a_s^{(j)}=0$ for all $j$, $1\le j\le p+q$. For $1\le j\le p$, let $m=(n_j-1)/4$ (resp., $m=(n_j+1)/4$) if $(n_j-1)/2$ is even (resp., $(n_j-1)/2$ is odd). Then $((\pi/2)+2m\pi)/n_j=\pi/2$ (resp., $((-\pi/2)+2m\pi)/n_j=\pi/2$). It follows from (4) that the interval $(\pi/2, (\pi/2)+(\pi/n_j))$ is contained in $S_j$ for all $j$. If $N=\max_{1\le j\le p}n_j$, then $(\pi/2, (\pi/2)+(\pi/N))\subseteq S_j$ for all $j$ and thus $\cap_{j=1}^pS_j\neq\emptyset$. We conclude from Corollary 4.7 that $d(A)=\sum_{j=1}^{p+q}d(A_j)+\sum_{k=1}^rd(B_k)$. The expression for $d(A)$ in terms of the $n_j$'s and $m_k$'s follows from Lemma 4.4 and Corollary 2.7.  $\blacksquare$

\vspace{5mm}

The final case we consider is for $A$ in Theorem 4.5 to have only two direct summands.

\vspace{5mm}

{\bf Proposition 4.9.} \emph{Let} $A=B\oplus C$, \emph{where}
$$B=\left[\begin{array}{cccc} 0 & b_1  & & \\ & 0 & \ddots & \\   & & \ddots & b_{m-1} \\ b_{m} & & & 0\end{array}\right] \ \ \ \ \ \mbox{\em and} \ \ \ \ \ C=\left[\begin{array}{cccc} 0 & c_1  & & \\ & 0 & \ddots & \\   & & \ddots & c_{n-1} \\ c_{n} & & & 0\end{array}\right]$$
\emph{with nonzero} $b_j$'\emph{s} \emph{and} $c_j$'\emph{s}. \emph{Then $d(A)=d(B)+d(C)+1$ if and only if $m=n$ is odd and $\sum_{j=1}^n(\arg b_j - \arg c_j)=(2\ell +1)\pi$ for some $\ell$}, $0\le\ell<n$. \emph{For the remaining case}, \emph{we have} $d(A)=d(B)+d(C)$.

\vspace{5mm}

{\em Proof}. Since $$d(B)+d(C)\le d(A)\le d(B)+d(C)+\lfloor \frac{1}{2}p\rfloor$$ by Corollary 4.7, where $p$ ($=0, 1$ or $2$) is the number of odd-sized matrices among $B$ and $C$, we obviously have $d(A)=d(B)+d(C)$ if either $m$ or $n$ is even. For the remaining part of the proof, we assume that both $m$ and $n$ are odd, and prove that (a) if $m\neq n$, then $d(A)=d(B)+d(C)$, and (b) if $m=n$, then $d(A)=d(B)+d(C)+1$ if and only if $|\beta-\gamma|=(2\ell+1)\pi$ for some $\ell$, $0\le\ell<n$, where $\beta=\sum_{j=1}^n\arg b_j$ and $\gamma=\sum_{j=1}^n\arg c_j$. As in (4), let
$$S=\left\{\begin{array}{ll}  \bigcup\limits_{\ell=-\infty}^{\infty}\big(\frac{1}{m}(\beta+\frac{\pi}{2}+2\ell\pi), \frac{1}{m}(\beta+\frac{3\pi}{2}+2\ell\pi)\big) & \mbox{     if } \, \frac{1}{2}(m-1) \mbox{  is even},\\
\bigcup\limits_{\ell=-\infty}^{\infty}\big(\frac{1}{m}(\beta-\frac{\pi}{2}+2\ell\pi), \frac{1}{m}(\beta+\frac{\pi}{2}+2\ell\pi)\big)& \mbox{   if } \, \frac{1}{2}(m-1) \mbox{   is odd}, \end{array}\right.$$
and let $T$ be defined analogously with $m$ and $\beta$ replaced by $n$ and $\gamma$, respectively.

\vspace{3mm}

To prove (a), note that $S'\equiv S\cap[0, 2\pi)$ and $T'\equiv T\cap [0, 2\pi)$ are such that $|S'|=|T'|=\pi$, and $|S'\cup T'|< 2\pi$ if $m\neq n$. Thus
$$|S'\cap T'|=|S'|+|T'|-|S'\cup T'|>\pi+\pi-2\pi=0.$$
and, therefore, $S\cap T\neq\emptyset$. Our assertion in (a) then follows from Corollary 4.7.

\vspace{3mm}

For the proof of (b), assume that $m=n$. In this case, it is easily seen that $S\cap T=\emptyset$ if and only if $|S\cap T|=0$, and the latter occurs exactly when $|\beta-\gamma|=(2\ell+1)\pi$ for some $\ell$, $0\le\ell<n$. Our assertion in (b) again follows from Corollary 4.7.  $\blacksquare$

\vspace{5mm}

\noindent
{\bf Acknowledgement}

We thank Chi-Kwong Li for some enlightening comments, which lead to a simplification of the proof of Theorem 3.3.

\end{document}